\documentclass[10pt]{article}

\usepackage{amsmath}
\usepackage{amsthm}
\usepackage{amsfonts}
\usepackage{amssymb}
\usepackage{amscd}
\usepackage{amsbsy}
\usepackage{graphicx}

\newcommand{\set}[1]{\,\left\{#1\right\}}
\newcommand{\setd}[2]{\,\left\{#1\ \colon\ #2\right\}}

\newtheorem{theorem}{Theorem}

\newtheorem{proposition}[theorem]{Proposition}
\newtheorem{definition}[theorem]{Definition}
\newtheorem{question}[theorem]{Question}
\newtheorem{remark}[theorem]{Remark}
\newtheorem{problem}[theorem]{Problem}

\newtheorem{conjecture}[theorem]{Conjecture}

\newcommand{\sign}{\operatorname{sign}}

\newcommand{\vertiii}[1]{{\left\vert\kern-0.25ex\left\vert\kern-0.25ex\left\vert #1 
    \right\vert\kern-0.25ex\right\vert\kern-0.25ex\right\vert}}

\newcommand{\RR}{\mathbb{R}}

\newcommand{\ZZ}{\mathbb{Z}}
\newcommand{\CC}{\mathbb{C}}
\newcommand{\NN}{\mathbb{N}}

\date{29 January 2014 }
\title{\ \\[0.4cm] \ \\ \bf  Group Actions on Banach Spaces}
\author{Piotr W. Nowak \footnote{Institute of Mathematics of the Polish Academy of Sciences, Warsaw, Poland  -- and -- University of Warsaw, Poland. 
Email: pnowak@mimuw.edu.pl. The author was partially supported by the Foundation for Polish Science}\hspace{2mm}}

\begin{document}

\maketitle


\thispagestyle{empty}

\begin{abstract}
\vskip 3mm\footnotesize{

\vskip 4.5mm
\noindent
We survey the recent developments concerning fixed point properties for group actions on Banach spaces. In the setting of Hilbert spaces 
such fixed point properties correspond to Kazhdan's property (T). Here we focus on the general, non-Hilbert case, we discuss 
the methods, examples and several applications.  
\vspace*{2mm}

\vspace*{2mm}
}

\end{abstract}
\vspace*{-8.2cm}

\vspace*{-17.6mm}\noindent{{\sl Handbook of\\
Group Actions}} \vskip8mm
\vspace*{7.6cm}

\section{Introduction}
Group actions are fundamental to understanding the geometry of both groups and spaces on which they act. 
Actions of groups on Banach spaces by affine endomorphisms,  that are additionally required to be  uniformly continuous, or isometric,
are particularly natural objects to study.

In addition to its geometric appeal, this topic has a natural connection with cohomology: various geometric properties of affine actions can be translated
into statements about the first cohomology group of $G$ with coefficients in the $G$-module formed by the Banach space $E$ with a representation 
$\pi$. This representation is the linear part of the affine action and cohomology group is denoted $H^1(G,\pi)$. Our main interest in this article will be the existence of fixed points for affine actions.
In the context of group cohomology, fixed point properties correspond to the vanishing of cocycles.

In the case when $\pi$ is a unitary representation on a Hilbert space the problem of establishing fixed point properties 
is motivated by property (T), introduced by Kazhdan. 
For a group $G$ property (T) is equivalent to the fact that every affine isometric action of $G$ 
on a Hilbert space has a fixed point. In other words, the cohomology
$H^1(G,\pi)$ vanishes for every unitary representation $\pi$ of $G$ on a Hilbert space.

Recently there has been growing interest in extending such rigidity properties to other Banach spaces, but even for such familiar classes 
as the Lebesgue spaces $L_p(\mu)$, or even spaces isomorphic to the Hilbert space, this program proved to be challenging. The lack of orthogonality 
presents a significant difficulty and new methods have
to be developed to prove fixed point properties in this general setting.

Our goal here is to give a fairly complete account of these recent developments and their applications. 
We purposely focus only on the case of Banach spaces which are not Hilbert
spaces, discussing the latter case mainly as motivation. In the case of Hilbert spaces Kazhdan's property (T) has been extensively studied and many excellent
sources are available, see for instance \cite{bekka-delaharpe-valette} and the references therein. 

We also discuss the opposite property of existence of a 
metrically proper affine isometric action on a Banach space, known in  the case of a Hilbert space as a-T-menability or the Haagerup property. 
An interesting phenomenon in this context is the existence of proper affine actions of hyperbolic
groups on $L_p$-spaces for $p\ge 2$ sufficiently large. The interplay between the existence of proper actions and the existence of fixed points
has applications to the geometry of  groups, for instance it allows to estimate various dimensions of boundaries of random hyperbolic groups.

We do not include proofs, instead we indicate, whenever possible, the ideas and methods behind the results. We also tried to include
a  comprehensive list of references.

\setcounter{tocdepth}{3}
\tableofcontents

\section{Preliminaries}

\subsection{Geometric properties of Banach spaces}
We will recall here several standard facts about Banach spaces and their geometry. This material can be found in 
many texts, we particularly recommend \cite{benyamini-lindenstrauss,johnson-lindenstrauss-handbook,megginson,wojtaszczyk}.

Let $V$ be a Banach space. By $V^*$ we denote the continuous dual of  $V$.
A particularly important class of examples is the class of Lebesgue spaces $L_p(\Omega, \mu)$ for $1\le p \le \infty$.
Many geometric features of Banach spaces can be expressed in terms of some measure of  convexity.

\begin{definition}
A  Banach space $(V, \Vert \cdot\Vert)$ is said to be strictly convex if
$$\left\Vert \dfrac{v+w}{2}\right\Vert < 1,$$
whenever $\Vert v\Vert=\Vert w\Vert=1$ and $v\neq w$.
\end{definition}
It is not hard to check that, for instance, the spaces $\ell_1(\Omega)$ and $\ell_{\infty}(\Omega)$, where $\Omega$ is a set with at least two elements,  
are not strictly convex. 
Similarly the space $c_0(\Omega)$ of functions vanishing at infinity, is not strictly convex. 
However, both $\ell_1(\Omega)$ and $c_0(\Omega)$ are separable and we have the following 

\begin{proposition}
Every separable Banach space admits an equivalent strictly convex norm.
\end{proposition}
Indeed, every separable Banach space $(V,\Vert\cdot\Vert_V)$ admits an injective operator $T:V\to H$ into a Hilbert space $H$ and the
norm $\Vert v\Vert'=\Vert v\Vert_V+\Vert Tv\Vert_H$ is a strictly convex norm on $V$, equivalent to $\Vert \cdot \Vert_V$.

When $1< p <\infty$ the Lebesgue spaces $L_p(\Omega, \mu)$ satisfy a stronger convexity property.

\begin{definition}
A  Banach space $(V, \Vert \cdot\Vert)$ is said to be uniformly convex if for every $\varepsilon>0$ there exists $\delta>0$ such that
$$\left\Vert \dfrac{v+w}{2}\right\Vert\le 1-\delta,$$
whenever $\Vert v\Vert=\Vert w\Vert=1$ and $\Vert v-w\Vert \ge \varepsilon$.
\end{definition}
The spaces $L_p(\Omega,\mu)$ for $1<p<\infty$ were shown to be uniformly convex by Clarkson \cite{clarkson}.
Another class of uniformly convex spaces is given by the Schatten $p$-class operators. A linear operator  $T:H_1\to H_2$, where $H_1, H_2$ are Hilbert spaces, 
is of Schatten class $p\ge 1$ if 
$$\operatorname{trace} \vert T\vert^p=\operatorname{trace}(T^*T)^{p/2}<\infty.$$ 
The space of such operators, denoted $\mathcal{C}_p$, is 
a Banach space with the norm $\Vert T\Vert _p=\left(\vert T\vert^p \right)^{1/p}$.
This definition can be generalized further to define non-commutative $L_p$-spaces, see e.g. \cite{pisier-xu}.

A Banach space, which admits an equivalent uniformly convex norm is called \emph{superreflexive}. Such spaces are automatically reflexive.

\begin{definition}
A Banach space $V$ is said to be uniformly smooth if 
$$\dfrac{\rho(t)}{t}\to 0,$$ 
as $t\to 0$, where 
$$\rho(t)=\sup\left\{\dfrac{\Vert v+w\Vert_V+\Vert v-w\Vert_V}{2}-1\ \Big\vert\ \Vert v\Vert_V=1, \Vert y\Vert \le t\right\},$$
is the modulus of smoothness of $V$.
\end{definition}
A Banach space $V$ is uniformly smooth if and only if its dual $V^*$ is uniformly convex. 
The three conditions of having an equivalent uniformly convex norm, an equivalent uniformly smooth norm, 
and an equivalent uniformly smooth and uniformly convex norm, are all equivalent.

A Banach space $V$ is said to have type $p\ge 0$, if there exists $C>0$ such that for any $n\in \NN$ and any collection of vectors
$v_1,\dots, v_n\in V$, the inequality
\begin{equation}\label{equation: type}
\dfrac{1}{\# A_n}\sum_{(\alpha_1,\dots,\alpha_n)\in A_n}\left\Vert\sum_{i=1}^n \alpha_iv_i\right\Vert_V^p\le C^p\sum_{i=1}^n\Vert v_i\Vert_V^p
\end{equation}
holds, 
where $A_n=\setd{(\alpha_1,\dots,\alpha_n)}{\alpha_i\in \set{-1,1}}$. There is a dual notion of cotype of a Banach space, and we refer to \cite[Chapter 11]{diestel-book}
for details. The inequality (\ref{equation: type}) is a weak version of the Hilbert space parallelogram law. Indeed, induction of the parallelogram law gives the equality
$$ \sum_{(\alpha_1,\dots,\alpha_n)\in A_n} \left\Vert\sum_{i=1}^n \alpha_iv_i\right\Vert^2=2^{n-1}\sum_{i=1}^n\Vert v_i\Vert^2,$$
for $A_n$ as above and any collection of vectors $\set{v_i}_{i=1}^{n}$ in the Hilbert space.
A superreflexive space has type $>1$. We refer to \cite[Chapter 11]{diestel-book} and \cite{johnson-lindenstrauss-handbook} 
for more details on type, cotype and their applications.

\begin{definition}
Let $V$ be a Banach space. A closed subspace $W\subseteq V$ is said to be complemented if there exists a closed subspace $W'\subseteq V$ such that
$V=W\oplus W'$.
\end{definition}

Equivalently, $W$ is complemented in $V$ if and only if there exists a bounded projection $P:V\to W$.
A theorem of Lindenstrauss and Tzafriri \cite{lindenstrauss-tzafriri} states that a Banach space, whose every closed subspace is complemented, 
is isomorphic to
the Hilbert space.

Given a family $\set{(V_i,\Vert \cdot\Vert_{V_i})}_{i\in I}$ of Banach spaces and $1\le p\le \infty$ we define the $p$-direct sum of the $V_i$, denoted
$$\left(\bigoplus_{i\in I} V_i\right)_p,$$ 
as the norm completion of the algebraic direct sum 
$\bigoplus_{i\in I} V_i$ in the norm 
$$\left\Vert \set{v_i}\right\Vert=\left(\sum_{i\in I} \Vert v_i\Vert_{V_i}^p\right)^{1/p}.$$ 
We have the identity 
$$\left(\bigoplus_{i\in I} V_i\right)_p^*=\left(\bigoplus_{i\in I} V_i^*\right)_q,$$
where $\dfrac{1}{p}+\dfrac{1}{q}=1$, see \cite{wojtaszczyk}.

A Banach space $V$ is said to be $L$-embedded, if $V$ is complemented in its second dual $V^{**}$ and $V^{**}=(V\oplus V')_1$, 
where $V'$ denotes the complement of $V$ in $V^{**}$.

\subsection{Representations, cocycles and cohomology}
Let $G$ be a discrete countable group and let $V$ be a Banach space. A \emph{representation of $G$ on $V$} is a homomorphism $\pi:G\to B_{inv}(V)$ into the group 
$B_{inv}(V)$ of bounded
invertible linear operators on
$V$.
The representation $\pi$ is said to be \emph{isometric} if $\pi_g$ is an isometry for each $g\in G$; $\pi$ is said to be  
\emph{uniformly bounded} if $\sup_{g\in G} \Vert \pi_g\Vert<\infty$.
By means of the representation $\pi$ we turn $V$ into a \emph{Banach $G$-module} and we can consider  the cohomology of $G$ with coefficients in the $G$-module $(V,\pi)$,
 usually referred to as \emph{cohomology with coefficients in $\pi$} and denoted $H^*(G,\pi)$. We refer to \cite{brown} for background on cohomology of groups.

We will focus on the 1-cohomology group $H^1(G,\pi)$, which is closely related to affine actions of $G$ on $V$. 
More precisely, an affine map $A:V\to V$ is a ``linear map+translation", 
$$Av=Tv+b,$$ 
for every $v\in V$, where $T$ is bounded linear operator on $V$ and $b\in V$. Given a representation $\pi$, an affine $\pi$-action of $G$ on $V$ 
is an action by affine maps  $A_g$ with the linear part
defined by the representation $\pi$; i.e., 
\begin{equation}\label{equation: affine action}
A_gv=\pi_g v+b_g,
\end{equation}
for every $v\in V$, where $b:G\to V$  satisfies the \emph{cocycle condition}, 
\begin{equation}\label{equation: cocycle condition}
b_{gh}=\pi_g b_h+b_g.
\end{equation}
Such a $b$ will be called a cocycle for $\pi$.
Another way to phrase this is that an affine action is a homomorphism into the semi-direct product $\operatorname{Aff}=V\rtimes B_{inv}(V)$.
We will be interested in properties of affine actions on various Banach spaces, with the linear part $\pi$ isometric or uniformly bounded.
A classical theorem of Mazur and Ulam \cite{mazur-ulam} states that a surjective metric isometry of a real Banach space is an affine isometry,
see \cite{nica-mazur-ulam} for the proof. A significant generalization of this result to certain ``fat'' groups of equi-continuous homeomorphisms 
was proved in \cite{mankiewicz}.

A cocycle $b:G\to V$ for $\pi$ is said to be a coboundary if there exists a vector $v\in V$, such that
$$b_g=\pi_gv-v,$$
for every $g\in G$. 
We denote by $Z^1(G,\pi)$ the linear space of cocycles for $\pi$, and the coboundaries form a space denoted $B^1(G,\pi)$.
\begin{definition}
The 1-cohomology group of $G$ with coefficients in $\pi$ is defined to be $H^1(G,\pi)=Z^1(G,\pi)\big/ B^1(G,\pi)$.
\end{definition}
It is easy to verify that $b$ is a coboundary satisfying $b_g=\pi_gv-v$ if and only if $v$ is a fixed point for the affine action determined by $\pi$ and $b$.
This simple observation leads to a geometric interpretation of the vanishing of 1-cohomology.
\begin{proposition}
$H^1(G,\pi)=0$ if and only if every affine action with linear part $\pi$ has a fixed point.
\end{proposition}
See also \cite{bekka-delaharpe-valette,nowak-yu} for discussions and background material.

\section{Fixed points}

\subsection{Kazhdan's property (T)}
Kazhdan's property (T) was introduced in \cite{kazhdan} and since then has been intensively studied. 
A unitary representation $\pi$ of a finitely generated group $G$ on a Hilbert space $H$ has almost invariant vectors if 
for every $\varepsilon>0$ there exists a non-zero vector $v\in H$ satisfying 
$$\Vert \pi_s v-v\Vert_H\le \varepsilon \Vert v\Vert_H,$$
for every generator $s\in G$.

\begin{definition}[Kazhdan]\label{definition: original definition of (T)}
A finitely generated group $G$ has property (T) if every unitary representation $\pi$ of $G$ that has almost fixed vectors
also has a non-zero fixed vector.
\end{definition}
It follows that an amenable group does not have property (T). Indeed, a characterization of amenability of a group $G$ is that there exists
a family $\set{f_i}_{i\in \NN}$, $f_i\in \ell_2(G)$, satisfying $\Vert f_i\Vert=1$ and $\Vert f_i-\lambda_s f_i\Vert\le 1/i$ for every generator $s\in G$, 
where $\lambda_s f(g)=f(s^{-1}g)$ is the left regular representation of $G$ on $\ell_2(G)$. 
(The $f_i$ can be taken to be the normalized characteristic functions of F\o lner sets). 
Thus the $f_i$ are a sequence of non-zero almost invariant vectors for $\lambda$, however it is clear that $\lambda$ does not have a non-zero invariant
vector.

In our setting it is much more natural to use a geometric  characterization of property (T), known as the the Delorme-Guichardet theorem (see \cite{bekka-delaharpe-valette}).

\begin{definition}\label{definition: (T) via fixed points}
A finitely generated group $G$ has Kazhdan's property (T) if and only if every affine isometric action of $G$ on a Hilbert space has a fixed point.
\end{definition}

As mentioned earlier, an affine action (\ref{equation: affine action}) with a cocycle $b$ has a fixed point $v$ if and only if $b_g=v-\pi_gv$ for all $g\in G$.
Consequently, according to the definition of 1-cohomology, we have the following reformulation.

\begin{proposition}
A finitely generated group $G$ has property (T) if and only if  $H^1(G,\pi)=0$ for every unitary representation $\pi$ of $G$.
\end{proposition}

Property (T) is a strong property, it is satisfied by relatively few groups but has numerous applications, including finite generation of lattices \cite{kazhdan}, solution of the
Ruziewicz problem \cite{rosenblatt,margulis,sullivan}, an explicit construction of expanders \cite{margulis-expanders}, various rigidity theorems for group actions and related operator algebras. We recommend 
\cite{bekka-delaharpe-valette,lubotzky} for a thorough introduction to property (T).

\subsection{Generalizing (T) to other Banach spaces}

We now arrive at the main object of interest in this survey.
Definition \ref{definition: (T) via fixed points} admits a natural generalization to normed spaces. Given a group $G$ and a Banach space $V$ we will thus 
be interested in the following property:

\begin{center}\label{fixed point property}
\emph{Every affine isometric action of $G$ on $V$ has a fixed point.}
\end{center}
Equivalently, 
\begin{center}
\emph{$H^1(G,\pi)=0$ for every isometric representation of $G$ on $V$.}
\end{center}

Even though the case when $V$ is the Hilbert space has been studied for several decades, 
for other Banach spaces such fixed point properties became an object of detailed study only recently. 
The articles \cite{fisher-margulis} and \cite{bfgm} were the first to systematically study fixed point properties on Banach spaces. 
One can immediately expect difficulties, as the proofs in the Hilbert space case rely heavily on the features of  inner product spaces such as
orthogonality, in particular existence of orthogonal complements,  and self-duality. In other Banach spaces neither of these is readily available. Every Banach space 
non-isomorphic to a Hilbert space contains a non-complemented closed subspace; 
computing dual spaces of closed subspaces of a Banach space $V$ means computing quotients  of
$V^*$; there is no standard notion of orthogonality between vectors in such vector spaces. All this turns statements which are evident 
in the Hilbert space setting 
into ones that are far from straightforward for other Banach spaces. 
In the next section we
recall a useful result, that also illustrates the
difficulties one has to deal with when working outside of Hilbert spaces. 

\subsection{Splitting off the invariant vectors} 
Observe first that in a Hilbert space $H$, the space of invariant vectors of a 
unitary representation $\pi$,  
$$H^\pi=\setd{v\in H}{\pi_gv=v\text{ for all } g\in G},$$
is closed and the representation $\pi$ preserves the orthogonal decomposition $H=H^\pi \oplus (H^\pi)^{\perp}$.

We would like to have a similar decomposition for any isometric representation of $G$ on any Banach space $V$.
However, this is not automatic: as mentioned earlier, such a $V$ always contains closed non-complemented subspaces whenever it is not isomorphic to the Hilbert space.

\begin{theorem}[\cite{bfgm}]
Let $\pi$ be an isometric representation of $G$ on a uniformly convex uniformly smooth Banach space $V$.
Then the subspace $V^{\pi}$ of invariant vectors is complemented in $V$.
\end{theorem}

The argument is based on the fact that for a uniformly smooth space $V$ every unit vector has a unique supporting functional $v^*$ of norm 1, which satisfies
$\langle v,w\rangle=\Vert v\Vert=1$, where $w\in V^*$, if and only if $w=v^*$. The (non-linear) 
map $S(V)\to S(V^*)$, $v\mapsto v^*$, between the unit spheres of $V$ and $V^*$,
intertwines a representation
$\pi$ on $V$ with its adjoint (contragradient) representation $\overline{\pi}$ on $V^*$, defined by 
$$\overline{\pi}_{\gamma}=\pi_{\gamma^{-1}}^*.$$
In particular, if $v\in V$ is fixed by $\pi$, then $v^*$ is fixed by $\overline{\pi}$. The complement of the space of $\pi$-fixed vectors in $V$ is then 
the annihilator of the $\overline{\pi}$-invariant vectors in $V^*$.

Recently the above result was generalized to a larger class of Banach spaces.
Given a representation $\pi$ of $G$ on a Banach space $X$ we will say that $\pi$ is a dual representation if there is a Banach 
space $Y$ and a representation $\rho$ of $G$ on $Y$ such that 
$$X=Y^* \ \ \ \ \ \text{ and } \ \ \ \ \ \pi=\overline{\rho}.$$
\begin{theorem}[\cite{garcia-nowak}]\label{theorem : garcia-nowak}
Let $\pi$ be a dual isometric representation on a dual Banach space $X$, such that $X^*$ is separable. Then 
$$X=X^{\pi} \oplus X_{\pi},$$
where  $X_{\pi}$ is a closed $\pi$-invariant subspace of $X$.
\end{theorem}
Theorem \ref{theorem : garcia-nowak} applies in particular to spaces $X$, which are duals of almost-reflexive Banach spaces, such as the well-known James space (see \cite{albiac-kalton}).

\subsection{Bounded orbits vs fixed points}\label{section: fixed points vs bounded orbits}

For a large class of Banach spaces the existence of a fixed point for an affine isometric action $A$ is equivalent to existence of a bounded orbit, namely existence of a vector $v$ 
satisfying
$$\sup_{\gamma\in G} \Vert  A_{\gamma} v\Vert<\infty.$$

If $V$ is a separable uniformly convex Banach space then any bounded non-empty set $K$ has a unique center, called the Chebyshev center. 
It is  the unique point $v\in V$ at which $\inf_{v\in V}\inf\setd{r>0}{K\subseteq B(v,r)}$ is attained. The existence of the Chebyshev center follows from the fact
that for reflexive spaces the weak topology and the weak$^*$ topology coincide, and, consequently, closed bounded sets are weakly compact. 
Uniqueness  is a consequence of uniform convexity. If now $K$ is the closed convex hull of a bounded orbit of an affine action then
the Chebyshev center of $K$ is a fixed point of the affine action.

If $V$ is a reflexive space then we can also apply the Ryll-Nardzewski to obtain a similar statement.
\begin{theorem}[The Ryll-Nardzewski fixed point theorem \cite{ryll-nardzewski}]
Let $V$ be a Banach space and $K\subset V$ be a non-empty weakly compact convex set. Then any group of affine isometries of $K$ 
has a fixed point.
\end{theorem}
In the setting of reflexive spaces the above theorem can be applied to the convex hull of a bounded orbit, producing a fixed point of the affine action.

Note that if an affine action of $G$ with linear part $\pi$ and a cocycle $b$ has two different fixed points $v\neq w$ in $V$, then 
$$\pi_g v-v=b_g=\pi_g w-w,$$
for all $g\in G$ and it follows that $v-w$ is a non-zero $\pi$-invariant vector. Consequently, if $\pi$ does not have non-zero
invariant vectors then fixed points are unique.

We remark that in general (i.e., for affine  actions on  non-reflexive spaces), the boundedness of an orbit does not imply that the corresponding cocycle is a coboundary.
The relation between the  boundedness of a cocycle and  being a coboundary for a  (uniformly bounded) representation $\pi$ 
is encoded in the bounded cohomology, $H^1_b(G,\pi)$.

\subsection{The fixed point spectrum}\label{section: the fixed points spectrum}

Let 
$$\mathcal{F}(G)=\left\{ p\in (1,\infty)\ \Big\vert
\begin{array}{l}
H^1(G,\pi)=0 \text{ for every isometric}\\
\text{representation } \pi \text{ on an } L_p\text{-space}
\end{array}\right\}$$
denote the \emph{fixed point spectrum of the group $G$}, with respect to the class of $L_p$-spaces. 
\begin{problem}
Given a finitely generated group $G$ with property (T), determine $\mathcal{F}(G)$.
\end{problem}
In the next section we will discuss the fact  that for a group with property (T), $(1,2]\subseteq \mathcal{F}(G)$.
As pointed out by C. Dru\c{t}u, in general it is only known that $\mathcal{F}(G)$ is open in $[1,\infty)$. 
When $p>2$ this can be proved using a similar argument as used in Proposition \ref{proposition: fisher-margulis}
below.
However, even the following natural 
question is open.

\begin{problem}[C.~Dru\c{t}u]
Is $\mathcal{F}(G)$ connected?
\end{problem}
For instance, if $G$ is a  hyperbolic group with property (T) then it is only known that there exists $\varepsilon=\varepsilon(G)>0$ such that 
$[2,2+\varepsilon)\subseteq \mathcal{F}(G)$ and that $ \mathcal{F}(G)$ is bounded, 
see Section \ref{subsection: hyperbolic groups act on L_p spaces}. In certain cases an estimate on $\varepsilon(G)$ can also be given,
see Section \ref{section: spectral criterion in the reflexive setting}.
\begin{problem}
Given $G$ with property (T), is there a critical value $p>2$ such that $\mathcal{F}(G)=(1,p)$?
\end{problem}
At present these problems are  far from being understood.

\subsection{The case $p\in [1,2)$}

We begin with a general result about the fixed point spectrum of groups with property (T).

\begin{theorem}[\cite{bfgm}]
Let $G$ have property (T). Then $(1,2]\subseteq \mathcal{F}(G)$.\\
Conversely, if $p\in (1,2)$ and $p\in \mathcal{F}(G)$ then $2\in\mathcal{F}(G)$.
\end{theorem}

To prove the first claim assume the contrary. The metric on $L_p$ for any $0<p\le 2$ is a negative definite function (see \cite{benyamini-lindenstrauss})
and the function $\varphi(\gamma)=\Vert b_{\gamma}\Vert_p$ is a negative definite function on $G$. However, a well-known 
characterization of property (T) states that any such $\varphi$ has to be bounded on a group with property (T).
The second claim follows from a result of Connes and Weiss \cite{connes-weiss}, see also \cite{glasner-weiss}.

The case $p=1$ is slightly different. The argument in the preceding paragraph implies that any affine isometric action of a property (T)
group has a bounded orbit, however, since $L_1$ is not reflexive unless it is finite-dimensional, we cannot deduce the existence of a fixed
point using the techniques described in Section \ref{section: fixed points vs bounded orbits}.
In fact, every group admits an isometric action without fixed points on a bounded convex subset of $L_1$.
An example of such an action is the translation action on the functions satisfying $\sum_{g\in G} f(g)=1$. 
In \cite{bader-gelander-monod} a fixed point theorem appropriate for this setting was proved.  
\begin{theorem}[\cite{bader-gelander-monod}]
Let $K\subseteq V$ be a non-empty, bounded subset of an $L$-embedded Banach spaces $V$. There exists a vector $v\in V$ such that
any isometry $\Phi$ of $K$ preserving $A$ satisfies $\Phi v=v$.
\end{theorem}
Since  $L_1$ is $L$-embedded, the above theorem covers the remaining case $p=1$.
Additional quantitative estimates and applications can be found in \cite{bader-gelander-monod}. 

\subsection{Unbounded spectrum}
There are several results which show that certain groups have fixed points for all isometric actions on $L_p$-spaces for all $1<p<\infty$; that is,
$\mathcal{F}(G)=[1,\infty)$. 
They are obtained by different methods, which we discuss below.

\subsubsection{The Mazur map and almost invariant vectors}\label{section: mazur map and higher rank groups}

The Mazur map $M_{p,q}:S(L_p(\Omega,\mu))\to S(L_q(\Omega,\mu))$ between the unit spheres of $L_p$ and $L_q$ is given by
$$M_{p,q}f(\omega)= \vert f(\omega)\vert^{\,p/q}\sign(f(\omega)),$$
for every $\omega\in \Omega$ and $f\in S(L_p(\mu))$.
The Mazur map is a uniform homeomorphism between the unit spheres of $L_p$, the modulus of continuity of $M_{p,q}$ and its inverse, $M^{-1}_{p,q}=M_{q,p}$ 
is given by the inequalities
$$\dfrac{p}{q} \Vert f-f'\Vert_{L_p} \le \Vert M_{p,q}f- M_{p,q}f'\Vert_{L_q}\le C \Vert f-f'\Vert_{L_p}^{p/q},$$
if $p<q$. 
See \cite{benyamini-lindenstrauss}.

As observed in \cite{bfgm}, the Mazur map can be used to produce almost invariant (respectively, invariant) vectors for an isometric representation on $L_p$ from 
almost invariant vectors (respectively, invariant) vectors of a representation on $L_q$.
It turns out that even tough the Mazur map is non-linear, it intertwines representations on $L_p$ spaces. The reason lies in the classification of isometries on 
$L_p$-spaces, proved  originally by Banach \cite{banach} and extended by Lamperti \cite{lamperti}. Given a (linear) isometry 
$U:L_p(\Omega, \mu)\to L_p(\Omega,\mu)$, 
where $1<p<\infty$ and $p\neq 2$, there exists a transformation $T:\Omega\to \Omega$, and a function $h:\Omega\to \CC$, such that 
$$(Uf)(x)=f(T(x)) h(x).$$
The Mazur map $M_{2,p}$ conjugates isometries of $L_p(\Omega,\mu)$ into isometries on $L_2(\Omega,\mu)$,
by the formula 
$$\pi'_g=M_{p,2}\,\pi_g\, M_{2,p}.$$ 
Moreover, it carries $\pi'$-invariant vectors to $\pi$-invariant vectors, and almost invariant vectors for $\pi$ 
to almost invariant  vectors for $\pi'$. This last fact follows from uniform continuity of $M_{p,q}$.
Consequently, we have
\begin{proposition}\label{proposition: (T) => almost fixed vectors definition on L_p}
Let $G$ be a group with property (T) and let $\pi$ be an isometric representation of $G$ on an $L_p$-space, $1<p<\infty$. If $\pi$ has
almost invariant vectors then $\pi$ has a non-zero invariant vector. 
\end{proposition}

The above is an $L_p$-version of the original definition of property (T) of Kazhdan (Definition \ref{definition: original definition of (T)}), for $2\le p<\infty$. 
However it turns out that
it is not enough to deduce vanishing of $H^1(G,\pi)$ for all isometric representations on $L_p$-spaces. As discussed later, hyperbolic groups with 
property (T) admit fixed point-free affine isometric actions on $L_p$-spaces for $p$ sufficiently large.

In the case of lattices in semisimple Lie groups additional properties allow to make the passage from the ``almost fixed vectors $\Rightarrow$ fixed vector" 
definition to vanishing of cohomology. One such property is a version of the Howe-Moore property for semisimple Lie groups and isometric 
representations on superreflexive Banach spaces, proved by Y.~Shalom (see the proof in \cite{bfgm}). 
Let now $k_i$, $i=1,\dots,m$ be local fields, $G_i$ be Zariski connected simple $k_i$-algebraic groups of $k_i$-rank at least $2$, and $G_i(k_i)$ be their $k_i$-points.
Denote  $G=\prod_{i=1}^m G_i(k_i)$. Such a $G$ will be called a higher rank group.

\begin{theorem}[\cite{bfgm}]
Let $\Gamma\subseteq G$ be a lattice in a higher rank group. Then $H^1(\Gamma,\pi)=0$ for every isometric representation on $L_p$, $1< p<\infty$.
\end{theorem}
In particular, the above theorem applies to $\operatorname{SL}_n(\ZZ)$, $n\ge 3$. A similar approach 
was used by M.~Mimura, who extended the above theorem to the so called universal lattices. 

\begin{theorem}[\cite{mimura}]
Let $k\in \NN$ and $n\ge 4$. Then $H^1\left(\operatorname{SL}_n(\ZZ[x_1,\dots,x_k]),\pi\right)=0$ for every isometric representation $\pi$
on $L_p$, $1< p<\infty$.
\end{theorem}

\subsubsection*{Non-commutative $L_p$-spaces}
The Mazur map method was also extended in a different direction. 
Puschnigg showed that for certain isometric representations $\rho$ of  a higher rank lattice $\Gamma$
on Schatten $p$-ideals $\mathcal{C}_p(H)$ on a Hilbert space $H$, 
$H^1(\Gamma,\rho)$ vanishes \cite{puschnigg}. The point is to apply the strategy outlined in \cite{bfgm}, 
together with a non-commutative version of the Mazur map.
Shalom's version of the Howe-Moore property also applies here, since $\mathcal{C}_p$ are uniformly convex.
The representations $\rho$ for which this argument goes through are the ones induced by a representation of $\Gamma$ on the Hilbert space $H$.

A further extension of this result was given by Olivier \cite{olivier}, who applied a non-commutative Mazur map to obtain the vanishing of  
$H^1(\Gamma,\rho)$ for every isometric representation of a higher rank lattice $\Gamma$ on a non-commutative $L_p$-space $L_p(\mathcal{M})$, where $\mathcal{M}$ is
a von Neumann algebra. A similar result due to Mimura for universal lattices $\operatorname{SL}_n(\ZZ[x_1,\dots,x_k])$, $n\ge 4$, can be found in \cite{mimura-pams}.

\subsubsection{Gromov monsters}
Gromov monsters are random groups, which, in a certain metric sense, contain expanders in their Cayley graphs, see \cite{arzhantseva-delzant,gromov-random}. 
They have exotic geometric features: they do not admit a coarse embedding into the Hilbert space. In \cite{naor-silberman} fixed point properties of such groups were
studied. 

Recall that for a graph $\Gamma$, the girth $g(\Gamma)$ denotes the length of the shortest non-trivial cycle.
Let $X$ be a metric space. We say that $X$ admits a sequence of high-girth $p$-expanders if there exists $k\in \NN$, constants $C,D>0$ and a sequence of $k$-regular
finite graphs $\set{\Gamma_n}$, where $\Gamma_n=(V_n,E_n)$, and $\vert V_n\vert \to \infty$, such that 
\begin{enumerate}
\item $g(\Gamma_n)\ge  C\log \vert V_n\vert$,
\item for every $f:V_n\to X$ the following Poincar\'{e} inequality holds:
$$\dfrac{1}{\vert V_n\vert^2}\sum_{v,w\in V_n} d_X(f(v),f(w))^p\le \frac{D}{\vert E_n\vert} \sum_{(v,w)\in E_n} d_X(f(v),f(w))^p.$$
\end{enumerate}
In other words, there exist $p$-expander graphs with respect to the geometry of $X$. 
In \cite{gromov-random} Gromov sketched a construction of groups which do not coarsely embed into the Hilbert space. 
These groups are constructed by introducing infinitely many random relations, which are modeled on relations in labeled expander graphs. 
We refer to \cite{arzhantseva-delzant,naor-silberman} for a detailed description of such groups.
A metric space $X$ is said to be $p$-uniformly convex if there exists $C>0$ such that for every triple of points $x,y,z\in X$, every geodesic segment $[yz]$,
and every $0\le t\le 1$ the inequality
$$d(x,[yz]_t)^p+Ct(1-t)d(y,z)^p\le (1-t)d(x,y)^p+td(x,z)^p$$
holds. Here $[yz]_t$ denotes the point on $[yz]$ at distance $td(y,z)$ from $y$.

\begin{theorem}[\cite{naor-silberman}]\label{theorem: naor-silberman}
Let $X$ be a $p$-uniformly convex metric space which admits a sequence of high girth $p$-expanders. Then for $d\ge 2$ and $k\ge 1$, 
with probability 1, every isometric action of a Gromov monster group on $X$ has a fixed point.
\end{theorem}
The proof rests on the idea that for a bounded orbit, a Chebyshev center is a fixed point of an action. If the orbit is unbounded one can try to
find certain averages over appropriate subsets of the group. 
It is proved in \cite{naor-silberman} that under the assumptions of the above theorem
such averages converge to a fixed point of the  action. Metric spaces to which Theorem \ref{theorem: naor-silberman} applies includes $L_p$-spaces, $1< p<\infty$, 
and
Euclidean buildings.

\subsubsection{Type $>1$}
A remarkable  result on the existence of fixed points on uniformly convex Banach spaces 
is due to V.~Lafforgue \cite{lafforgue-jta}, who introduced in \cite{lafforgue-duke} 
a version of strengthened property (T) 
and studied its various applications.
For a group $G$, property (T) can be characterized by the existence of a self-adjoint idempotent $p$ in the maximal $C^*$-algebra $C^*_{max}(G)$,
such that for every unitary representation $\pi$ of $G$ on a Hilbert space $H$, the image of $\pi(p)$ consists of the vectors fixed by $\pi$.

Consider now representations $\pi:G\to B(H)$, satisfying
\begin{equation}\label{equation: lafforgues condition on representations}
\Vert \pi_\gamma\Vert \le e^{l(\gamma)},
\end{equation}
where $l( \gamma)$ denotes  the length of $\gamma\in G$.
Equip the complex group ring $\CC G$ with the norm
$$\Vert f\Vert_l= \sup \Vert \pi(f)\Vert_{B(H)},$$
where the supremum is taken over all representations $\pi$ satisfying (\ref{equation: lafforgues condition on representations}).
The algebra $C_{l}^*(G)$ is defined to be the completion of $\CC G$ under the norm $\Vert\cdot\Vert_l$.
The algebra $C_l^*(G)$ is a Banach *-algebra.
Observe, that in this language, $C_{max}^*(G)$ is the algebra associated to the trivial length function.

\begin{definition}[\cite{lafforgue-duke}]
A group $G$ is said to have property (TR) if for every length function $l$ on $G$ there exists a constant $s>0$ such that for every $c>0$
there exists a self-adjoint idempotent $p\in C_{sl+c}(G)$, such that  for every 
representation $\pi$ of $G$ on a Hilbert space $H$, satisfying $\Vert \pi(\gamma)\Vert\le e^{sl(\gamma)+c}$, 
the image of $\pi(p)$ consists of the vectors fixed by $\pi$.
\end{definition}

There are several interesting facts about this notion, among them it was shown that infinite hyperbolic groups do not 
have property (TR). On the other hand, uniform lattices in $\operatorname{SL}_3(\RR)$ do have property (TR).

More generally, Lafforgue  also introduces a strengthened version of property (TR) with respect to the class of Banach spaces of type $>1$. 
He shows that the group $\operatorname{SL}_3(\mathbb{F})$ and its uniform lattices, where $\mathbb{F}$ is a non-archimedean local field, have 
such a strong property (TR)
with respect to Banach spaces with  type $>1$. In a subsequent article \cite{lafforgue-jta} the following fixed point property was proved.
\begin{theorem}
Let $G$ be a group with strong property (TR) with respect to the class of Banach spaces of type $>1$. 
Then $H^1(G,\pi)=0$ for every  isometric representation of $G$ on a Banach space with type $>1$.
\end{theorem}
In particular, if $\mathbb{F}$ is a non-archimedean local field then for a uniform lattice $G\subseteq \operatorname{SL}_3(\mathbb{F})$, 
$H^1(G,\pi)=0$ for every isometric 
representation of $G$ on a Banach space with type $>1$. Lafforgue's methods were used recently in \cite{liao} 
to extend these results to 
simple algebraic groups of higher rank.

Recently de la Salle \cite{dls} showed that $\operatorname{SL}_3(\RR)$ has Lafforgue's strong property (T) with respect to a large class of Banach spaces, 
defined using type and cotype (see \cite{dls} for precise definitions). 
Subsequently de Laat and de la Salle \cite{dlsdl} generalized this result to simple Lie groups of higher rank.

\subsection{Bounded spectrum}
A group $G$ for which $\mathcal{F}(G)=[1,\infty)$ can be thought of as a very rigid group. However, from the point of view of 
the program outlined in Section \ref{section: the fixed points spectrum}, the case $\mathcal{F}(G)\neq [1,\infty)$ seems to be much more difficult.
Nevertheless certain methods are available and we discuss them below.

\subsubsection{A general result}
We begin with a general statement  by D.~Fisher and G.~Margulis. The argument was included in \cite{bfgm}.
\begin{proposition}\label{proposition: fisher-margulis}
Let $G$ be a group with property (T). Then there exists $\varepsilon(G)>0$, such that $H^1(G,\pi)=0$ for any isometric representation $\pi$ 
of $G$ on  $X=L_p(\mu)$ whenever $2\le p<2+\varepsilon(G)$.
\end{proposition}
The proof is based on taking an ultralimit of  the spaces $L_{p}(\mu)$, and of the corresponding actions, as $p\to 2$, together with certain quantitative estimates.
As observed by C.~Dru\c{t}u, the same argument also proves that $\mathcal{F}(G)$ is an open subset of $[1,\infty)$.

\subsubsection{The spectral criterion in the reflexive setting}\label{section: spectral criterion in the reflexive setting}

The spectral method for proving property (T) is often referred to as the geometric method. The main idea,
originating from the work of Garland \cite{garland}, is that if a group $G$ is acting on a simplicial 2-dimensional complex, in which 
links of all vertices have first positive eigenvalue of the
Laplacian strictly greater than $1/2$, then $G$ has property (T). For actions on Hilbert spaces this method was used in 
\cite{ballmann-swiatkowski,pansu,zuk1,zuk2}. A more general version, relying on the energy of harmonic maps, applies to metric spaces in general and 
has been used in \cite{wang-1,wang-2} and discussed in \cite[3.11]{gromov-random}.

In \cite{nowak-poincare} a version of the Garland method was extended from Hilbert spaces to reflexive Banach spaces. To state the theorem we first discuss
Poincar\'{e} inequalities and link graphs.

Let $\Gamma=(\mathcal{V},\mathcal{E})$ be a finite graph, $V$ a Banach space and $1\le p \le \infty$. The $p$-Poincar\'{e} inequality is the inequality
\begin{equation}\label{equation: Poincare inequality}
\sum_{v\in \mathcal{V}} \Vert f(v)-Mf\Vert_V^p\le \kappa \sum_{v\in \mathcal{V}}\sum_{w\in \mathcal{V}, v\sim w} \Vert f(v) -f(w)\Vert_V^p,
\end{equation}
where $Mf=\dfrac{1}{\#\mathcal{V}} \sum_{v\in \mathcal{V}}f(v)\deg(v)$ is the mean value of $f:\mathcal{V}\to V$.

\begin{definition}
The Poincar\'{e} constant $\kappa(p,V,\Gamma)$ of the graph $\Gamma$ is the optimal constant $\kappa$ for which the inequality (\ref{equation: Poincare inequality})
holds for all $f:\mathcal{V}\to V$.
\end{definition}
For $p=2$ and $V=\RR$ the Poincar\'{e} constant is the square root of the inverse of $\lambda_1$, the first non-zero eigenvalue of 
the discrete Laplacian on $\Gamma$. Thus the Poincar\'{e} constant is a natural generalization of the spectral gap to the non-Hilbertian setting. 
Also note that by integration, $\kappa(p,\RR,\Gamma)=\kappa(p,L_p(\Omega,\mu),\Gamma)$, for any measure space $(\Omega,\mu)$. Estimating Poincar\'{e} constants 
for various $p$ and $V$ is a difficult task even for relatively simple graphs.

\begin{definition}
Let $G$ be a group, generated by a finite generating set $S=S^{-1}$, not containing the identity element. 
The link graph $\mathcal{L}(S)$ is a finite graph, defined as follows:
\begin{enumerate}
\item vertices of $\mathcal{L}(S)$ are elements of $S$,
\item two vertices $s,t\in S$ are connected by an edge if $s^{-1}t\in S$ (equivalently, $t^{-1}s\in S$).
\end{enumerate}
\end{definition}
The link graph is, in general, a non-regular graph, and for an arbitrary generating set it may turn out to be disconnected.
However, given any generating set $S$, the link graph of the generating set $S'=(S\cup S^2)\setminus \set{e}$ is connected.
The eigenvalues of the Laplacian of the link graph have significant influence on the cohomology of the group.
Given $p$ we define $p^*$ by the relation $1/p+1/p^*=1$.

\begin{figure}
\begin{center}
\includegraphics[width=100pt]{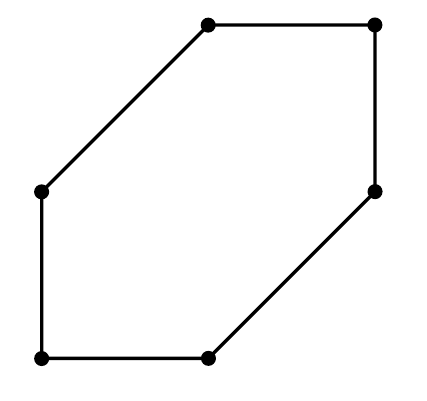}
\caption{The link graph $\mathcal{L}(S)$ of the generating set $S=\set{(1,0),(-1,0),(0,1),(0,-1),(1,1),(-1,-1)}$ in $\mathbb{Z}^2$}
\end{center}
\end{figure}

\begin{theorem}[\cite{nowak-poincare}]\label{theorem: nowak spectral criterion reflexive spaces}
Let $G$ be a group, generated by a finite, symmetric set $S$, not containing the identity element
and let $V$ be a reflexive Banach space.
If  
$$\kappa(p,V,\mathcal{L}(S))< 2^{1/p}\ \ \ \ \ \ \text{ and } \ \ \ \ \ \ \kappa(p^*,V^*,\mathcal{L}(S))\le 2^{1/p^*}, $$
then $H^1(G,\pi)=0$ for every isometric representation $\pi$ of $G$ on $V$. 
\end{theorem}

Note that when $p=p^*=2$ and $V$ is a Hilbert space, then we obtain precisely the spectral criterion in the form proved by \.{Z}uk \cite{zuk2} (see also 
\cite{nowak-yu} for a proof of \.{Z}uk's result).
We remark that although the method presented below and the one used to prove  Theorem  \ref{theorem: naor-silberman} both rely on Poincar\'{e} inequalities,
they are different in nature: the latter applies to very exotic groups, but is insensitive to the exact value of the Poincar\'{e} constant, while the method of 
Theorem \ref{theorem: nowak spectral criterion reflexive spaces}
described below requires an estimate of the Poincar\'{e} constant, but applies to well-behaved groups, such as hyperbolic groups.

The proof of the above theorem is entirely linear in nature and allows to obtain some additional information about the behavior of cocycles
for isometric representations on reflexive Banach spaces. It also allows to estimate Kazhdan-type constants
for representations on $L_p$-spaces, and consequently eigenvalues of the discrete $p$-Laplacian on finite quotients of groups.

Roughly speaking, the proof relies on methods of continuous homological algebra. If $\pi$ is an isometric representation of $G$
on a Banach space $V$, then showing that every cocycle is a coboundary amounts to showing 
surjectivity of the continuous operator $dv(g)=\pi_gv-v$ from the $V$ onto the space of cocycles $B^1(G,\pi)$.
By duality, this corresponds to the adjoint operator $d^*:B^1(G,\pi)^*\to E^*$ being bounded below; that is, $d^*$ is injective and has closed image.
The calculations in \cite{nowak-poincare} are performed on the restrictions of cocycles to the generating set $S$. In that setting
the explicit formula for $d^*$ can be computed and it turns out that  $d^*f$ can be expressed in terms of  $Mf$, the mean value of $f$ over the generators, for appropriate $f$. 
This allows to build a connection with Poincar\'{e} inequalities on the link graph. 

Most importantly however, theorem \ref{theorem: nowak spectral criterion reflexive spaces}
allows to obtain the first explicit lower bounds on $\varepsilon(G)$ such that $[2,2+\varepsilon(G))\subset \mathcal{F}(G)$
for certain groups, which we now discuss

\subsubsection*{$\widetilde{A}_2$ groups} 
$\widetilde{A}_2$-groups are groups which act transitively on the vertices of a building of type $\widetilde{A}_2$ \cite{cmsz}.
They are parametrized by powers $q=k^n$ for  prime $k$ and we denote them by $G_q$.
The fact that these groups have property (T) was first proved in \cite{cms}, however it also follows from the spectral method.
Indeed, these groups admit natural presentations for which the link graph is the incidence graph of a finite projective plane and these graphs have sufficiently 
large spectral gaps.
Certain $\widetilde{A}_2$ groups are lattices in higher rank groups and the results of \cite{bfgm}, described in Section \ref{section: mazur map and higher rank groups}, apply,
giving $\mathcal{F}(G)=[1,\infty)$ for these groups. However, for other $\widetilde{A}_2$-groups the only known method of proving fixed point properties
for $L_p$-spaces is the 
spectral method of Theorem \ref{theorem: nowak spectral criterion reflexive spaces}.

The linear Poincar\'{e} constants $\kappa(p,V,\Gamma)$  allow for interpolation and applying Theorem \ref{theorem: nowak spectral criterion reflexive spaces} we obtain
\begin{theorem}[\cite{nowak-poincare}]\label{theorem : fixed points for A_2-groups}
Let $q=k^n$ for some $n\in \NN$ and prime number $k$. Then $$H^1(G_q,\pi)=0$$
for any isometric representation $\pi$ of $G_q$ on $L_p(\Omega,\mu)$ 
for every
$$2\le p<\dfrac{\ln(q^2+q+1)+\ln(q+1)}{\dfrac{1}{2}\ln(2(q^2+q+1)(q+1))-\ln(2)-\ln\left(\sqrt{1-\dfrac{\sqrt{q}}{q+1}}\right)}.$$
\end{theorem}

Although we can only prove that a certain  explicit, bounded interval is contained
in $\mathcal{F}(G)$, it is likely that the following question has an affirmative answer.
\begin{question}
Let $G$ be an $\widetilde{A}_2$-group. Is it true that $\mathcal{F}(G)=[1,\infty)$?
\end{question}

\subsubsection*{Hyperbolic groups} The Gromov density model for random groups allows to produce examples of hyperbolic groups with property (T). 
A group $G$ in Gromov's density model $\Gamma(n,k,d)$ is defined by a presentation $G=\langle S \vert R\rangle$, where the generating set $S$ satisfies 
$\# S=n$, and the set of relations $R$ consists of $(2n-1)^{kd}$ relations, chosen independently and uniformly from the set of all relations of length $k$. Given $n\in \NN$ 
and the density $0\le d\le 1$, we say that a property $P$ holds almost surely for groups in $\Gamma(n,k,d)$ if 
$$\lim_{k\to\infty} \mathbb{P}(G\in \Gamma(n,k,d)\text{ satisfies } P)=1.$$

Gromov proved that if  $d<1/2$ then the group we obtain is almost surely hyperbolic.
On the other hand, if $d>1/3$ then the resulting random group almost surely has property (T)  \cite{zuk2}, see also \cite{kotowscy}. 
Applying interpolation one obtains an explicit number $C_p(G)>0$, for which $[2,2+C_p(G))\subseteq \mathcal{F}(G)$. In the case of hyperbolic groups this has an additional consequence.
\begin{theorem}[\cite{nowak-poincare}]\label{theorem: conformal dimension estimate}
Let $G$ be a random hyperbolic group in the Gromov density model  with $1/3< d<1/2$. Then $\operatorname{confdim}(\partial G)\ge C_p(G)$.
\end{theorem}
The problem of estimating the conformal dimension of the boundary for random hyperbolic groups was posed by Gromov in \cite{gromov-asymptotic}.
We discuss the conformal dimension briefly in section \ref{subsection: hyperbolic groups act on L_p spaces}.
Theorem \ref{theorem: conformal dimension estimate} 
is a consequence of Theorem \ref{theorem: nowak spectral criterion reflexive spaces} and a result due to Bourdon and Pajot \cite{bourdon-pajot}
that hyperbolic groups act without
fixed points on $L_p$-spaces, for $p\ge \operatorname{confdim}(\partial G)$, see section \ref{subsection: hyperbolic groups act on L_p spaces} for a discussion.

The only other result in the direction of estimates of the conformal dimension of the boundaries of random hyperbolic groups is due to J.~Mackay \cite{mackay},
however it is very different in spirit. The random groups considered there are constructed using  
densities $d<1/16$ and are a-T-menable by a result of Ollivier and Wise \cite{olliver-wise} (see Section \ref{section: metrically proper actions}).

\begin{remark}[Poincar\'{e} constants.]\normalfont
An important step in obtaining lower estimates of the right hand end of the fixed point spectrum in the above method is an estimate of the Poincar\'{e} constant.
In general, various interpolation techniques can be applied here. The numerical values for $\varepsilon(G)$ obtained in \cite{nowak-poincare} for
$\widetilde{A}_2$-groups and random hyperbolic groups are between $2$ and $2.2$. A very broad open problem is the computation the optimal 
Poincar\'{e} constants  for various graphs, different values of $p$ and certain Banach spaces $V$. Specific examples of $V$ for which there are virtually no
estimates of Poincar\'{e} constants include Schatten $p$-class operators $\mathcal{C}_p$ and $p$-direct sums of $\ell_q$, $p\neq q$.
\end{remark}

We also would like to mention that a non-linear approach to Garland-type results, outlined by Gromov in \cite[3.11]{gromov-random}, was also used in the context of 
Banach spaces in \cite{bourdon,oppenheim-spectral}.

\subsection{Uniformly bounded representations}
Extending cohomological vanishing  to isometric representations on Banach spaces other than the Hilbert space is one direction of generalizing property (T),
another is the extension of  fixed point properties from isometric to a larger class of affine actions on the Hilbert space. Particularly 
interesting is the class of uniformly bounded representations. A representation $\pi:G\to B(H)$ into the invertible operators on 
some Hilbert space is said to be uniformly bounded if 
$$\sup_{g\in G}\Vert \pi_g\Vert<\infty.$$
A simple renorming of $H$,
\begin{equation}\label{equation: renorming uniformly bounded}
\vertiii{v}=\sup_{g\in G}\Vert \pi_gv\Vert,
\end{equation}
allows to view uniformly bounded representations on $H$ as isometric representations on Banach spaces isomorphic to $H$.
In this case the following  question is natural.
\begin{question}
Let $G$ be a group with property (T). Is it true that $H^1(G,\pi)=0$ for every uniformly bounded representation on a Hilbert space?
\end{question}

In the case of lattices in higher rank groups and universal lattices the answer is affirmative \cite{bfgm,mimura}. 
Also in the case of Gromov monsters the argument in \cite{naor-silberman} can be used to show that such a fixed point property holds.
However, in the case of other groups with property (T) things are much less clear. An unpublished result of Y.~Shalom states that 
the group $\operatorname{Sp}(n,1)$
has a uniformly bounded representation $\pi$ on a Hilbert space, for which  $H^1(G,\pi)\neq 0$.  
The only other result in this direction concerns random hyperbolic groups.

\begin{theorem}[\cite{nowak-poincare}]
Let $G$ be a random group in the Gromov density model with $1/3<d<1/2$. Then  almost surely $H^1(G,\pi)=0$ for every uniformly bounded 
representation on the Hilbert space, satisfying
$$\sup_{g\in G}\Vert \pi_g\Vert<\sqrt{2}.$$
\end{theorem}
The above is an application of Theorem \ref{theorem: nowak spectral criterion reflexive spaces} applied to the Banach space $H$ with the 
equivalent norm (\ref{equation: renorming uniformly bounded}) and $p=2$. A similar statement can be proved for $\widetilde{A}_2$-groups.
Theorem \ref{theorem: nowak spectral criterion reflexive spaces} also applies to uniformly bounded representations on $L_p$-spaces for certain $p$.
The following problem is open.
\begin{conjecture}[Y.~Shalom]
Let $G$ be a non-elementary hyperbolic group. There exists a uniformly bounded representation of $G$ on a Hilbert space,
for which $H^1(G,\pi)\neq 0$ and for which there exists a proper cocycle.
\end{conjecture}

Recently the spectral condition for vanishing of $L_2$-cohomology was extended to uniformly bounded representations in
\cite{koivisto}, generalizing \cite{nowak-poincare}. It allows to use the spectral criterion to show the vanishing of cohomology of groups acting on simplicial complexes
and can be applied to the examples of groups constructed in \cite{ballmann-swiatkowski}.

\subsection{Reduced cohomology}
To finish the discussion of fixed points we would like to briefly mention the reduced cohomology and a growing body of work in the direction of 
vanishing of the reduced cohomology with values in uniformly bounded representations on reflexive Banach spaces. 

Let $\pi$ be a uniformly bounded 
representation of a group $G$ on a Banach space. The space of $\pi$-cocycles, $Z^1(G,\pi)$, can be equipped with the topology
of uniform convergence on compact subsets: a net of cocycles $b_\alpha$ converges to $b$ if for every compact subset $K\subset G$
the restriction $(b_\alpha)_{\vert K}$ converges to $b_{\vert K}$ in norm. (In the  setting of finitely generated groups it suffices to consider convergence in 
norm on the generators.) We denote by $\overline{B^1(G,\pi)}$ the closure of $B^1(G,\pi)$.
\begin{definition}
The reduced cohomology $\overline{H}^1(G,\pi)$ is  defined to be the quotient $Z^1(G,\pi)\big/\overline{B^1(G,\pi)}$.
\end{definition}
Observe that  for a group generated by a finite set $S$ the vanishing 
of $\overline{H}^1(G,\pi)$ is equivalent to the fact that for every $\pi$-cocycle $b$, there exists a sequence of vectors 
$v_n\in V$ such that
$$\Vert b_s-(\pi_sv_n-v_n)\Vert\to 0,$$
for every generator $s\in G$.
Rearranging the terms we obtain
$$\Vert (\pi_sv_n+b_s)-v_n\Vert_V=\Vert A_sv_n-v_n\Vert\to 0,$$
where $A_g v=\pi_g v+b_\gamma$ is the affine actions associated to $\pi$ and $b$. In other words,
the sequence $\set{v_n}$ forms a sequence of  \emph{almost fixed points} of the affine action $A$.

Vanishing of the reduced cohomology was studied by Shalom in the setting of unitary representations on Hilbert spaces \cite{shalom}.
For $V=\ell_p(G)$ the problem of vanishing the reduced cohomology is closely related to the problem of vanishing of the reduced $\ell_p$-cohomology.
\begin{conjecture}[Gromov]
Let $G$ be an amenable group. Then the reduced $\ell_p$-cohomology of $G$ vanishes for $1<p<\infty$.
\end{conjecture}
The conjecture is motivated by the case $p=2$: the vanishing of the reduced $\ell_2$-cohomology of amenable groups is a classical results of Cheeger and
Gromov \cite{cheeger-gromov}.
We refer to \cite{gromov-asymptotic} for an overview of this topic and to \cite{bourdon-martin-valette,martin-valette} for results on the case $p\neq 2$.

More recently some  results on reduced cohomology have been extended to the setting of reflexive spaces, see \cite{bader-rosendal-sauer,nowak-lpcoho}.
A fundamental technique in these considerations has become the Ryll-Nardzewski fixed point theorem \cite{ryll-nardzewski}.
 
\section{Metrically proper actions}\label{section: metrically proper actions}

\subsection{a-T-menability}

a-T-menability was defined  by Gromov, see \cite{gromov-asymptotic}, and independently in a different setting by  Haagerup \cite{haagerup}.
It can be defined in terms of properties of affine isometric actions.

\begin{definition}
Let $\pi$ be a representation of $G$ on a Banach space $V$.
An affine isometric action with linear part $\pi$ is metrically proper if $\lim_{g\to \infty}\Vert b_g\Vert_V\to \infty$.
\end{definition}
Clearly, if $G$ admits a metrically proper affine isometric action with linear part $\pi$ then the cohomology $H^1(G,\pi)$ cannot vanish, since
the coboundaries $\pi_gv-v$ are always bounded cocycles.

\begin{definition}
A finitely generated group is called a-T-menable (or is said to have the Haagerup property) if it admits a proper affine isometric action on a Hilbert space.
\end{definition}

Examples of a-T-menable groups include amenable groups \cite{bekka-cherix-valette} and free groups \cite{haagerup}, 
a-T-menability is also preserved by free products. We refer to \cite{cherix-et-al,nowak-yu} and the references therein for more details.
Aside from the relation to property (T),  significant interest in a-T-menability is a consequence of a remarkable result of Higson and Kasparov.
\begin{theorem}[\cite{higson-kasparov}]
Let $G$ be a finitely generated group. If $G$ is a-T-menable then the Baum-Connes conjecture holds for $G$.
\end{theorem}

The Baum-Connes conjecture, whenever it is true, is a vast generalization of the Atiyah-Singer index theorem and has several important consequences,
including the Novikov conjecture. See \cite[Chapter 8]{nowak-yu} and \cite{valette-bcc}. a-T-menability is also an important  property from the point
of view of large scale geometry, we refer to \cite{nowak-yu} for details.

\subsection{The case $p\in[1,2)$}
We again begin with a general result.
Similarly as for existence of fixed points, the behavior for $L_p$-spaces between $1\le p< 2$ is similar to the case $p=2$.
\begin{theorem}[\cite{nowak, nowak-corrected}]
A finitely generated group $G$ is a-T-menable,
if and only if $G$ admits a proper affine isometric action on $L_p[0,1]$ for every $1\le p \le 2$.

If $G$ is a-T-menable then $G$ admits a proper affine isometric action on $L_p[0,1]$ for every $1\le p<\infty$.
\end{theorem} 
The proof is an application of negative definite functions and a dynamical characterization of a-T-menability due to P.~Jolissaint \cite{cherix-et-al}.
A similar fact for $\ell_p$-spaces was proved in  \cite{cdh} using the geometry of median spaces and actions on spaces with measured walls.

\subsection{Amenable groups}
A very versatile method of constructing proper cocycles was presented by M.~Bekka, P.-A.~Cherix and A.~Valette 
in \cite{bekka-cherix-valette}. They answered a question of Gromov, proving the following
\begin{theorem}[\cite{bekka-cherix-valette}]
An amenable group $G$ admits a proper, affine isometric action on a Hilbert space (i.e., is a-T-menable).
\end{theorem}
The proof can be easily adapted to $\ell_p$-spaces and in fact shows that any amenable group admits a proper affine isometric action on $\ell_p(X)$ for a countable
set $X$ and for any $1\le p<\infty$. F\o lner's characterization of amenability states that a group $G$ is amenable if and only if there exists 
a sequence of finite set $F_n\subseteq G$ such that 
$$\dfrac{\# F_n\triangle sF_n}{\# F_n}\le\dfrac{1}{2^n},$$
for every generator $s\in G$. 
By defining $f_n=\dfrac{1_{F_n}}{\Vert 1_{F_n}\Vert_p}$ we construct the isometric action on $V=\left(\bigoplus_{n\in \NN}\ell_p(G)\right)_p$
as follows. The representation on $V$ is $\pi_g=\bigoplus_{n\in \NN}\lambda_g$, where $\lambda_g$ is the left regular representation of $G$ on $\ell_p(G)$:
$$\lambda_gf(x)=f(g^{-1}x),$$
for every $g,x\in G$.
The cocycle is then defined by the formula
$$b_g=\bigoplus_{n\in \NN} f_n-\lambda_gf_n.$$
One now verifies that $b$ is well-defined and proper, both of these facts follow from the properties of the F\o lner sequence $\set{F_n}$. Since $V$ is isometrically 
isomorphic to $\ell_p(X)$ for a countable set $X$, the theorem follows.

\subsection{Reflexive spaces}
In the same spirit as with property (T), we can generalize a-T-menability to any Banach space $V$, by asking if a group $G$ 
admits a metrically proper affine isometric action on $V$.
It turns out that already reflexivity is broad enough to allow a proper isometric action of any finitely generated group.

\begin{theorem}[\cite{brown-guentner}]
Any finitely generated group $G$ acts properly by affine isometries on the Banach space $V=\left(\bigoplus_{p=2,3,\dots} \ell_p(G)\right)_2$.
\end{theorem}
The idea for the proof is to show that every discrete group allows a version of amenability with respect to the Banach space $c_0(G)$: for every 
$n\in \NN$ the function $f(g)=\max\set{ 1-\dfrac{d(e,g)}{n},0}$, where $e$ denotes the identity elements, is finitely supported, has $c_0$-norm 1 and satisfies
$$\Vert f-\lambda_s\cdot f\Vert_{c_0}\le \dfrac{1}{n},$$
for every generator $s\in G$. Finite support allows to approximate the $c_0$-norm of $\Vert f-\lambda_g f\Vert_{c_0(G)}$ by $\Vert f-\lambda_g\cdot f\Vert_{\ell_p(G)}$ 
as $p\to \infty$. Adapting the previous argument one then shows that
the formula $b_g=\bigoplus f-\lambda_g f$ defines a proper cocycle for the isometric representation $\bigoplus \lambda$ on $\left(\bigoplus_{p\in \NN} \ell_p(G)\right)_2$.

\subsection{Hyperbolic groups and $L_p$-spaces}\label{subsection: hyperbolic groups act on L_p spaces}

A particularly  interesting case is that of hyperbolic groups. The simplest example of a hyperbolic group, the free group $\mathbb{F}_n$ on
$n$ generators, is a-T-menable. In fact, one can construct an example of a proper affine isometric action using the geometric properties of the Cayley graph of $\mathbb{F}_n$.

Consider the space $\ell_2(E_{\pm})$, where  $E_{\pm}$ denotes the set of pairs $(s,t)$, where $s^{-1}t$ or $t^{-1}s$ is a generator of $\mathbb{F}_n$. In other words,
$E_{\pm}$ is the space of oriented edges of the Cayley graph associated to a symmetric set of generators. The group acts naturally on $E_{\pm}$, inducing 
a representation of $\mathbb{F}_n$ on $\ell_2(E_{\pm})$. We define a  cocycle $b_g$ for this representation to be the characteristic function of the 
union of all oriented edges on the unique path connecting $g$ to the origin. It can be easily verified that $b$ is a proper cocycle.  Another argument
to prove a-T-menability of $\mathbb{F}_2$ uses the GNS construction. It is easy to prove that the square of the metric on the tree, i.e., a Cayley graph of $\mathbb{F}_2$, is 
a proper negative definite function. Via the GNS construction such functions induce proper affine isometric action on a Hilbert space.

The same strategy fails for general hyperbolic groups. Indeed, as discussed earlier, 
some hyperbolic groups are known to have Kazhdan's property (T) and every affine isometric action on a Hilbert space has bounded orbits.  
However, as it turns out,
every hyperbolic group $G$ admits a proper isometric action on an $L_p$-spaces, for $p=p(G)>2$ sufficiently large.

The first result, that certain hyperbolic groups have fixed point free actions on $L_p$-spaces for certain sufficiently large $p>2$, follows from the work Pansu \cite{pansu-95}, 
who showed that
the $L_p$-cohomology of the group $\operatorname{Sp}(n,1)$ does not vanish in degree 1 for $p\ge 4n+2$. It can be easily seen (see e.g., \cite{martin-valette}) that this implies
the existence of a fixed point free affine action (i.e., a non-trivial cocycle) associated with the regular representation on $L_p(G)$ for $p\ge 4n+2$. 
Note that since $\operatorname{Sp}(n,1)$ has property (T), the above  also implies that for $p\neq 2$ the two generalizations of property (T),
are not equivalent: the fixed point property is stronger than the obvious generalization of Definition \ref{definition: original definition of (T)}. 

Generalizing the case of $\operatorname{Sp}(n,1)$, M.~Bourdon and H.~Pajot proved  \cite{bourdon-pajot} that every non-elementary hyperbolic group has non-vanishing $\ell_p$-cohomology for
$p\ge \operatorname{confdim}\partial G$, the conformal dimension of $G$. Recall that the conformal dimension of the boundary of a hyperbolic group was defined by 
P.~Pansu \cite{pansu-conformal} to be the number 
$$\operatorname{confdim}\partial G=\inf\left\{
\dim_H(\partial G,d)\ \Big\vert
\begin{array}{l} d \text{ is quasi-conformally}\\ 
\text{equivalent to a visual metric}\end{array}\right\},$$
where $\dim_H$ denotes the Hausdorff dimension. We refer to \cite{kapovich-benakli} for an overview of the conformal dimension in the context of boundaries of groups.

The non-existence of fixed points can be strengthened to existence of a proper cocycle. 
The first result  was proved by G.~Yu.

\begin{theorem}[\cite{yu-hyperbolic}]\label{theorem: yu hyperbolic proper action}
Let $G$ be a hyperbolic group. There exists $p\ge 2$, which depends on $G$, such that $G$ acts properly by affine isometries on 
$\ell_p(G\times G)$.
\end{theorem}

A version of the above theorem in the special case of fundamental groups of hyperbolic manifolds (the ``classical" hyperbolic groups) can be found in \cite{nowak-yu}.
The proof of Theorem \ref{theorem: yu hyperbolic proper action} in \cite{yu-hyperbolic}
relies on an averaging construction of I.~Mineyev \cite{mineyev-gafa,mineyev}. This construction allows to find discrete analogs 
of tangent vectors. 

A new construction has been given recently by B.~Nica. 
\begin{theorem}[\cite{nica}]\label{theorem: hyperbolic bogdan}
Let $G$ be a non-elementary hyperbolic group. Then $G$ acts properly by isometries on the space $L_p(\partial G\times \partial G)$,
for every $p\ge \hbar(G)$, the hyperbolic dimension of $G$.
\end{theorem}

Actions of a hyperbolic groups on their boundaries are well-behaved, even if the group itself is not (e.g., has property (T)). 
Nica's beautiful construction also relies on Mineyev's averaging \cite{mineyev}, however in a different way than Yu's proof. 
Namely, he uses a new class of visual metrics on the boundary, that were constructed in \cite{mineyev}, and that behave better in
certain aspects than the standard visual metrics induced by the word length metric.

One of the main points is a construction of an analog of the 
Bowen-Margulis measure on the product $\partial G\times\partial G$. This measure is infinite but $G$-invariant. 
The  construction of the cocycle is geometric in nature. The hyperbolic dimension 
of $\Gamma$, $\hbar(\Gamma)$, introduced in \cite{mineyev},
is modeled on the conformal dimension of the boundary, discussed earlier.

\begin{question}
Does there exist a non-hyperbolic group $G$ with property (T), such that $\mathcal{F}(G)$ is bounded?
\end{question}

\section{Final remarks}
\subsection{Groups of homeomorphisms}
We end with two open questions, both of which concern groups of homeomorphisms of compact manifolds.

\begin{question}(G.~Yu)
Let $G$ be a finitely generated subgroup of the diffeomorphism group of a closed manifold $M$. Does $G$ admit a proper affine isometric action on 
a uniformly convex Banach space?
\end{question}
Even in the case of $M=S^1$ this problem is extremely interesting. For instance in \cite{navas} (see also \cite{bekka-delaharpe-valette}) it was proved that
among groups of sufficiently smooth diffeomorphisms of $S^1$ there are no infinite Kazhdan groups.
\begin{theorem}\cite{navas,navas2}
Let $G$ be a group with $[2,2+c(G))\subseteq \mathcal{F}(G)$. For any $\alpha>\dfrac{1}{c(G)}$ and any 
homomorphism $\varphi:G\to \operatorname{Diff}^{1+\alpha}_+(S^1)$, the image $\varphi(G)$ is a finite cyclic group.
\end{theorem}

Another question concerns mapping class groups. It is natural to ask if they admit proper affine isometric actions on uniformly
convex Banach spaces. In fact we have the following attractive conjecture.
\begin{conjecture}[B.~Nica]
The mapping class group of a surface admits a proper affine isometric action on an $L_p$-space for some sufficiently large $p=p(G)\in [2,\infty)$.
\end{conjecture}

\subsection{Relative fixed point properties}

Let $G$ be a group and $Q\subseteq G$ be a non-compact subset. For a Banach space $V$ 
we can define a fixed point property relative to $Q$ by requiring that
any affine isometric action $A$ of $G$ on a Banach space $V$ has a point fixed by $Q$:
$$A_gv=v,$$
for every $g\in Q$. It is customary to consider $Q=H$ to be a subgroup of $G$. For Hilbert space this property is known as the relative property (T).
For other Banach spaces such  properties have not been studied systematically.

\subsection*{Acknowledgments}
I would like to thank Uri Bader, Bogdan Nica and Rufus Willett for comments and suggestions which greatly improved the exposition.

\end{document}